\documentclass[12pt,reqno]{amsart}

\usepackage{amsthm,amsmath,amsfonts,amssymb,graphicx,graphics}
\usepackage{color}

\newtheorem{thm}{Theorem}
\newtheorem{lem}[thm]{Lemma}
\newtheorem{prop}[thm]{Proposition}
\newtheorem{cor}[thm]{Corollary}
\theoremstyle{definition}

\def\Sym{\mathfrak{S}}
\def\Der{{D}}
\def\Desar{{K}}
\def\DEZ{\mathop{\rm DEZ}\nolimits}
\def\IDES{\mathop{\rm IDES}\nolimits}
\def\DES{\mathop{\rm DES}\nolimits}
\def\des{\mathop{\rm des}\nolimits}
\def\maj{\mathop{\rm maj}\nolimits}
\def\exc{\mathop{\rm exc}\nolimits}
\def\iexc{\mathop{\rm iexc}\nolimits}
\def\fix{\mathop{\rm fix}\nolimits}
\def\pix{\mathop{\rm pix}\nolimits}

\def\lec{\mathop{\rm lec}\nolimits}

\def\inv{\mathop{\rm inv}\nolimits}

\title[Permutations with Extremal number of Fixed Points]
{Permutations with Extremal number of Fixed Points}

\author[Guo-Niu Han \and Guoce Xin]{Guo-Niu Han$^1$ \and Guoce Xin$^2$}
\address{$^{1,2}$ Center for Combinatorics, LPMC-TJKLC, Nankai University,
Tianjin, 300071, P. R. China}
\address{$^1$ I.R.M.A. UMR 7501, Universit\'e Louis Pasteur et
CNRS,  7 rue Ren\'e-Descartes, F-67084 Strasbourg, France}

\email{$^1$guoniu@math.u-strasbg.fr, $^2$gxin@nankai.edu.cn}
\date{June 22, 2007}
\begin{document}

\begin{abstract}
We extend Stanley's work on alternating permutations with extremal
number of fixed points in two directions: first, alternating
permutations are replaced by permutations with a prescribed descent
set; second, instead of simply counting permutations we study their
generating polynomials by number of excedances. Several techniques
are used: D\'esarm\'enien's desarrangement combinatorics, Gessel's
hook-factorization and the analytical properties of two new
permutation statistics ``DEZ" and ``lec". Explicit formulas for the
maximal case are derived by using 
symmetric function tools.
\end{abstract}

\maketitle \thispagestyle{empty}




\section{Introduction} 
Let $J=\{j_1, j_2, \ldots, j_r\}_<$ be a set of integers arranged
increasingly and let $\Sym_J$ denote the set of all permutations on $J$.
For each permutation $\sigma=\sigma(j_1)\sigma(j_2)\cdots \sigma(j_r) \in
\Sym_J$ define the \emph{number of excedances}, the \emph{number of
fixed points} and the \emph{descent set} of $\sigma$ to be
\begin{align}
\fix \sigma&=|\{i: 1\leq i\leq r, \sigma(j_i)=j_i\}|, \nonumber \\
 \exc \sigma&=|\{i:
1\leq i\leq r, \sigma(j_i)>j_i\}|, \nonumber \\
 \DES \sigma&=\{i: 1\leq i\leq r-1,
\sigma(j_i)>\sigma(j_{i+1}) \}, \label{e-descent}
\end{align}
respectively.
A permutation without fixed point is called a {\it derangement}.
When $J=[n]:=\{1, 2, \dots , n\}$, we recover the
classical definitions. The set $\Sym_{[n]}$ is abbreviated by $\Sym_n$, and
for $\sigma\in \Sym_n$ we write $\sigma_i$ for $\sigma(i)$.
Our main results are the following Theorems~\ref{t-1-gen}
and \ref{exc-gen}.

\begin{thm}\label{t-1-gen}
Let $J$ be a subset of $[n-1]$.

\text{\rm (i)} If $\sigma\in\Sym_n$ and $\DES\sigma=J$, then
$$
\fix\sigma\leq n-|J|.
$$

\text{\rm (ii)} Let $F_n(J)$ be the set of all permutations $\sigma$
of order $n$ such that $\DES\sigma=J$ and $\fix\sigma=n-|J|$.
Furthermore, let $G(J)$ be the set of all derangements $\tau$ on $J$
such that $\tau(i)>\tau(i+1) $ whenever $i$ and $i+1$ belong to $J$.
Then
$$
\sum_{\sigma\in F_n(J)} s^{\exc\sigma} =\sum_{\tau\in G(J)}
s^{\exc\tau}.
$$
\end{thm}

\medskip
\noindent {\bf Example.} For $n=8$ and $J=\{1,2,3,6\}$, there are
two permutations in $F_n(J)$, both having two
excedances: 74315628 and 74325618.
On the other hand, there are two derangements in $G(J)$, both
having two excedances: 6321 and 6312.

\begin{thm}\label{exc-gen}
Let $\Der_0^J(n)$ be the set of all derangements $\sigma$ on $[n]$
such that $\DES\sigma=J$, and let $\Der_1^J(n)$ be the set of all
permutations $\sigma
$
on $[n]$ such that $\DES \sigma
=J$ with exactly one fixed
point. If $J$ is a proper subset of $[n-1]$, then there is a
polynomial $Q_n^J(s)$ with positive integral coefficients such that
$$
\sum_{\sigma\in\Der_0^J(n)} s^{\exc\sigma}
-\sum_{\sigma
\in\Der_1^J(n)} s^{\exc\sigma
} =(s-1) Q_n^J(s).
$$
\end{thm}

\medskip
\noindent {\bf Example.} For $n=6$ and $J=\{1,3,4,5\}$ there are six
derangements in $\Der_0^J(n)$:
$$
216543, 316542, 416532, 436521, 546321, 645321;
$$
and there are six permutations in $\Der_1^J(n)$:
$$
326541, 426531, 516432, 536421, 615432, 635421.
$$
The numbers of excedances are respectively 3, 3, 3, 4, 3, 3 and 3,
3, 2, 3, 2, 3, so that
$$
\sum_{\sigma\in\Der_0^J(n)} s^{\exc\sigma}
-\sum_{\sigma
\in\Der_1^J(n)} s^{\exc\sigma
} =(5s^3+s^4) -
(4s^3+2s^2)= (s-1)(s^3+2s^2).
$$

Theorem \ref{t-1-gen} extends Stanley's work on alternating
permutations (that we explain next) with maximal number of fixed
points, and Theorem \ref{exc-gen} extends the corresponding minimal
case. The extensions are in two directions:
first, alternating permutations
are replaced by permutations with a prescribed descent set; second, instead
of simply counting permutations we study their
generating polynomials by number of excedances.

A permutation $\pi=\pi_1\pi_2\cdots \pi_n\in \Sym_n$ is said to be
{\it alternating} (resp. {\it reverse alternating}) if $\pi_1 >
\pi_2 < \pi_3 > \pi_4 < \dots $ (resp. if $\pi_1 < \pi_2 > \pi_3 <
\pi_4 > \dots $);
or equivalently, if $\DES \pi$ is $\{1,3,5,\dots\} \cap
[n-1]$ (resp. $\{2,4,6,\dots\} \cap [n-1] $). Therefore, results on
permutations with a prescribed descent set apply to alternating
permutations. Let $\Der_k(n)$ be the set of permutations in $\Sym_n$
with exactly $k$ fixed points. Then $\Der_0(n)$ is the set of
derangements of order~$n$. Write $d_k(n)$ (resp. $d^*_k(n)$) for the
number of alternating (resp. reverse alternating) permutations in
$\Der_k(n)$.
The next two corollaries are immediate consequences of
Theorems 1 and 2.

\begin{cor}[\cite{St}, Conjecture 6.3]\label{t-1}
Let $D_n$ denote the number of derangements.  Then, for $n\geq 2$ we have
$$
d_n(2n)
=d_{n+1}(2n+1)
=d^*_{n+1}(2n+1)
=d^*_{n+2}(2n+2)=D_{n}.
$$
\end{cor}

\begin{cor}[\cite{St}, Corollary 6.2]
\label{t-2} For $n\geq 2$ we have $d_0(n)=d_1(n)$ and
$d_0^*(n)=d_1^*(n)$.
\end{cor}

Stanley enumerated $D_k(n)$ and came up with Corollaries \ref{t-1}
and \ref{t-2} on alternating permutations with extremal number of
fixed points. He then asked for combinatorial proofs of them. This
is the motivation of the paper. The results in Corollary \ref{t-1},
conjectured by Stanley, was recently proved by Chapman and Williams
\cite{Wi} in two ways, one directly and the other using the newly
developed concept of permutation tableaux \cite{SW}. In Section
\ref{sec-cor3} we give a direct proof of a generalized form
of Corollary \ref{t-1}. Corollary \ref{t-2} is actually a special
case of a more general result due to Gessel and Reutenauer, which
itself can be derived from Theorem 2 by setting $s=1$, as stated in
the next corollary.

\begin{cor}[\cite{gessel-reutenauer}, Theorem 8.3]\label{t-2-gen}
Let $J$ be a proper subset of $[n-1]$. Then, the number of
derangements in $\Sym_n$ with descent set $J$ is equal to the number
of permutations in $\Sym_n$ with exactly one fixed point and descent
set $J$.
\end{cor}

The paper is organized as follows. In Section 2 we give the proof of
Theorem \ref{t-1-gen} that contains the results for the maximal case.
Section 3 includes a direct proof of an extension of Corollary \ref{t-1}. 
Section 4 introduces the necessary part of Gessel and Reutenauer's work for
enumerating the maximal case. Section 5 is devoted to the proof of
Theorem \ref{exc-gen} dealing with  the minimal case.
We conclude the paper by making several remarks of analytic nature (see
Section 6). In particular, Corollary \ref{exc-sum-n-q1}, proved combinatorially, should deserve
an analytic proof.
Several techniques are used:
D\'esarm\'enien's desarrangement combinatorics \cite{De84},
Gessel's hook-factorization \cite{Ge91}
and the analytical properties of
two new permutation statistics ``DEZ" and ``lec"
\cite{FH06, FH07}.


\section{Permutations with maximal number of fixed points}\label{sec-th1} 
Our task in this section is to prove Theorem 1. The proof relies on
the properties of the new statistic ``$\DEZ$" introduced by Foata and Han \cite{FH06}.
For a permutation $\sigma=\sigma_1\sigma_2\cdots \sigma_n\in
\Sym_n$ let~$\sigma^0=\sigma^0_1 \sigma^0_2 \cdots\sigma^0_n$
be the word derived from $\sigma$ by
replacing each fixed point $\sigma_i=i$ by~$0$. The set-valued
statistic ``$\DEZ$" is defined by
$$
\DEZ \sigma =\DES \sigma^0:=\{i: 1\leq i\leq n-1, \sigma^0_i>\sigma^0_{i+1} \}.
$$
For example, if $\sigma=
8\,2\,1\,3\,5\,6\,4\,9\,7$, then $\DES\sigma=\{1,2,6,8\}$,
$\sigma^0=8\,0\,1\,3\,0\,0\,4\,9\,7$ and
$\DEZ\sigma=\DES\sigma^0=\{1,4,8\}$. The basic property of the
statistic ``DEZ" is given in the following proposition.

\begin{prop}[\cite{FH06}, Theorem 1.4]\label{DEZequi}
The two three-variable statistics $(\fix, \exc, \DEZ)$ and $(\fix,
\exc, \DES)$ are equi-distributed on the symmetric group $\Sym_n$.
\end{prop}
More precisely, Proposition \ref{DEZequi} asserts that there is a
bijection $\Phi: \Sym_n \mapsto \Sym_n$ such that
$$\fix \pi =\fix \Phi(\pi), \quad \exc \pi=\exc \Phi(\pi), \quad \DES \pi =\DEZ \Phi(\pi),  \text{ for all } \pi\in \Sym_n. $$
By Proposition \ref{DEZequi} Theorem \ref{t-1-gen} is equivalent to
the following Theorem \ref{t-1-gen'}, where the statistic ``DES" has
been replaced by ``DEZ".

\addtocounter{thm}{-1}
\renewcommand\thethm{\ref{t-1-gen}$'$}
\begin{thm}\label{t-1-gen'}
Let $J$ be a subset of $[n-1]$.

\text{\rm (i)} If $\sigma\in\Sym_n$ and $\DEZ\sigma=J$, then
$$
\fix\sigma\leq n-|J|.
$$

\text{\rm (ii)} Let $F'_n(J)$ be the set of all permutations
$\sigma$ of order $n$ such that $\DEZ\sigma=J$ and
$\fix\sigma=n-|J|$. Furthermore, let $G(J)$ be the set of all
derangements $\tau$ on $J$ such that $\tau(i)>\tau(i+1) $ whenever
$i$ and $i+1$ belong to $J$. Then
$$
\sum_{\sigma\in F'_n(J)} s^{\exc\sigma} =\sum_{\tau\in G(J)}
s^{\exc\tau}.
$$
\end{thm}
\renewcommand\thethm{\arabic{thm}}

\begin{proof}[Proof of Theorem \ref{t-1-gen'}]
Let $\sigma$ be a permutation such that $\DEZ\sigma=J$ and let $i\in J$. Then
$\sigma^0_i>\sigma^0_{i+1}\ge 0$, so that $i$ is not a fixed
point of $\sigma$. It follows that $\sigma$ has at least $|J|$ non-fixed
points. This proves (i).

Now, consider the case where $\sigma$ has exactly $n-|J|$ fixed points.
Then $J$ is the set of all the non-fixed points of $\sigma$. By removing the
fixed points from $\sigma$ we obtain a derangement~$\tau$ on $J$. If
$i,i+1\in J$, then $\tau(i)=\sigma(i)>\sigma(i+1)=\tau(i+1)$.
It follows that $\tau\in G(J)$. On the other hand,
take any derangement $\tau \in G(J)$ and let $\sigma$ be the permutation
defined by
$$\sigma(i)=
\begin{cases}
\tau(i),   &\text{if $i\in J$,}\cr
i,         &\text{if $i\not\in J$.}\cr
\end{cases}
$$
Then $\DEZ\sigma=J$. It is easy to see that $\sigma\in F'_n(J)$ and
$\exc\sigma=\exc\tau$.
This proves the second part of Theorem \ref{t-1-gen'}.
\end{proof}

\medskip
\noindent {\bf Example.} Suppose $n=8$ and $J=\{1,2,3,6\}$. Let us
search for the permutations $\sigma\in \Sym_8$ such that $\fix\sigma=8-|J|=4$ and
$\DEZ\sigma=J$. There are two derangements $\tau$ in $G(J)$, namely, $6321$
and $6312$, both having two excedances, so that the two corresponding
elements $\sigma$ in $F'_n(J)$ are $63245178$ and $63145278$, both
having two excedances.

\medskip
\noindent\textbf{Remarks.} (i) For permutations with descent set $J$ it
is easy to show that the maximum number of fixed points is $n-|J|$,
except when $J$ consists of an odd number of consecutive integers.
In the latter exceptional case the only decreasing permutation has exactly
one fixed point and therefore is not a derangement.

(ii) The first part of Theorem 1 can also be proved directly
by using
the fact that in any consecutive decreasing
subsequence of $\pi$, say $\pi_i>\pi_{i+1}>\cdots
>\pi_{i+k}$, there is at most one fixed point in $\{i,i+1,\dots,i+k\}$.
However the ``DEZ" statistic is an essential tool in the proof of the
second part.


\section{An extension of Corollary 3}\label{sec-cor3} 
Stanley's conjectured result
in Corollary \ref{t-1} was first proved by Williams \cite{Wi} using
the newly developed concept of permutation tableaux. A direct proof
without using permutation tableaux was later included in her updated
version with Chapman. Our direct proof 
was independently derived just
after Williams' first proof. It has the advantage of automatically
showing the following extension (Proposition \ref{p-conj-ext}).
We only give the generalized form for $d_n(2n)=D_n$, since the
other cases are similar. All of the three proofs are bijective, and
the bijections are all equivalent.
Note that Proposition \ref{p-conj-ext} is still a corollary of
Theorem \ref{t-1-gen}.

\begin{prop}\label{p-conj-ext}
The number of alternating permutations in $\Sym_{2n}$ with $n$ fixed
points and $k$ excedances is equal to the number of derangements in
$\Sym_n$ with $k$ excedances.
\end{prop}

Let $\pi$ be an alternating permutation. Then, each doubleton
$\{\pi_{2i-1}, \pi_{2i}\}$ contains at most one fixed point. This
proves the following lemma.

\begin{lem}\label{l-group}
Each alternating permutation $\pi \in \mathfrak{S}_n$ has at most
$\lceil n/2 \rceil$ fixed points. When this maximum is reached,
either $2i-1$, or $2i$ is a fixed point of $\pi$
$(2\le 2i \le n+1)$.
\end{lem}

When the underlying set of the permutation $\pi$ is not
necessarily $[n]$, we use $\pi(i)$
instead of $\pi_i$ for convenience. An integer $i$ is called an
\emph{excedance}, a \emph{fixed point}, or a \emph{subcedance} of
$\pi$ if $\pi(i)>i$, $\pi(i)=i$, or $\pi(i)<i$, respectively.
\begin{proof}[Proof of Proposition \ref{p-conj-ext}]
Let $\pi \in \Sym_{2n}$ be alternating and have exactly $n$ fixed
points. It follows from Lemma \ref{l-group} that for each $i$ we
have the following property: either $2i-1$ is a fixed point and $2i$
a subcedance, or $2i-1$ is an excedance and $2i$ a fixed point.
Conversely, if the property holds, the permutation $\pi$ is
necessarily alternating, because $\pi(2i)\le 2i <2i+1\le \pi(2i+1)$;
$\pi(2i-1)\ge 2i-1\ge \pi(2i)-1$. Those inequalities imply that
$\pi(2i-1)>\pi(2i)$, since $2i-1$ and $2i$ cannot be both fixed
points.

By removing all fixed points of $\pi$ we obtain a derangement
$\sigma$ on an $n$-subset of $[2n]$. The standardization of
$\sigma$, which consists of replacing the $i$-th smallest element of $\sigma$ by
$i$, yields a derangement $\tau$ on $[n]$. We claim that the map
$\varphi:\pi \mapsto \tau$ is the desired bijection.
Since the standardization preserves excedances, subcedances and
fixed points, it maps one element of $\{\pi(2i-1), \pi(2i)\}$  to
$\tau(i)$. It follows that $\tau(i)>i$ if and only if $2i-1$ is
an excedance and $2i$ is a fixed point of $\pi$, and that $\tau(i)<i$
if and only if $2i-1$ is a fixed point and $2i$ is a subcedance of
$\pi$. Thus, the set of all fixed points of $\pi$ can be constructed
from $\tau$. The map $\varphi$ is then reversible.

The proposition then follows since the bijection preserves the
number of excedances.
\end{proof}

\medskip
\noindent {\bf Example.} Let
$\pi=\begin{pmatrix}
1 & 2        & 3 & 4 & 5 &            6 & 7 & 8 & 9 & 10 \\
3  & \bar{2} & 6 & \bar{4} & \bar{5} & 1 & 10 & \bar{8} & \bar{9} & 7 \\
\end{pmatrix}$.
Removing all the fixed points gives
$\sigma=\begin{pmatrix}
1 &  3 &   6 & 7 &   10 \\
3  &  6 &   1 & 10   & 7 \\
\end{pmatrix},$
standardized to
$\tau=\begin{pmatrix}
1 &  2 &   3 & 4 &   5 \\
2  &  3 &   1 & 5   & 4 \\
\end{pmatrix}.$
Conversely, $\tau$ has excedances at positions $1,2,4$ and subcedances
at positions $3,5$.
This implies that $2,4,8 $ and $5,9$ are fixed points of $\pi$ and
hence we can construct $\pi$. 
Furthermore, we have $\exc\pi=\exc\sigma=\exc\tau=3$.


\section{Enumeration for the maximal case} 

In this section we will use Theorem \ref{t-1-gen} to enumerate the
number of permutations with a prescribed descent set and having the
maximal number of fixed points.
Every descent set $J\subseteq [n-1]$ can be partitioned into blocks of
consecutive
integers, such that numbers from different blocks differ by at least
$2$. Let $J^b=(a_1,a_2,\dots,a_k)$ denote the sequence of the size of
the blocks.
For instance, if $J=\{1,2,3,6\}$, then $1,2,3$ form a block and $6$
itself forms  another block. Hence $J^b=(3,1)$.
Let $M_J$ denote the number of derangements in
$\Sym_n$ with descent set $J$ having $n-|J|$ fixed points.
By Theorem~\ref{t-1-gen} the number $M_J$ depends only on $J^b$.
Thus, we can denote $M_J$ by $M(a_1,a_2,\dots, a_k)$.
\begin{thm}
\label{t-max-fix}
The number $M(a_1,\dots,a_k)$ is the coefficient of
$x_1^{a_1}\cdots x_k^{a_k}$ in the expansion~of
$$
\frac{1}{(1+x_1)(1+x_2)\cdots (1+x_k)(1-x_1-x_2-\cdots -x_k)}.
$$
\end{thm}
An immediate consequence of Theorem \ref{t-max-fix} is the following
Corollary \ref{c-d-symmetry}, which says that $M(a_1,a_2,\dots,a_k)$
is symmetric in the $a_i$'s.
\begin{cor}\label{c-d-symmetry}
For each permutation $\tau \in \Sym_k$ we have
$$M(a_1,a_2,\dots,a_k)=M(a_{\tau_1},a_{\tau_2},\dots, a_{\tau_k}). $$
\end{cor}
For example, $M(3,1)$ counts two derangements $4312$ and $4321$;
$M(1,3)$ counts two derangements $3421$ and $4321$. This symmetry
seems not easy to prove directly.
Using Theorem \ref{t-max-fix} an
explicit formula for $M(a_1, a_2,\ldots, a_k)$
can be obtained when
$k=1, 2$.
We have
$M(a)=1$ if $a$ is even, and $M(a)=0$ if $a$ is
odd; also
$$ M(a,b)=\sum_{j=2}^{a+b} \sum_{i=0}^j (-1)^j \binom{a+b-j}{a-i}.$$
\medskip

To prove Theorem \ref{t-max-fix} we need some notions
from \cite{gessel-reutenauer}, where Gessel and Reutenauer
represented the number of permutations with given cycle structure
and descent set by the scalar product of two special characters of
the symmetric group introduced by Foulkes
\cite{foulkes-up-down,foulkes-eulerian}. Their results were also key
ingredients in \cite{St} for the enumeration of alternating
permutations by number of fixed points.
In what follows, we assume the basic knowledge of symmetric
functions (see, e.g., \cite{Lasc,Macd,EC2}).
The scalar product $\langle \ , \ \rangle$
of two symmetric functions
is a bilinear form
 defined for all partitions $\lambda$ and $\mu$
by
\begin{equation}
\label{e-scalar} \langle m_\lambda, h_\mu \rangle =\langle h_\mu,
m_\lambda \rangle= \delta_{\lambda\mu},
\end{equation}
where $m_\lambda$ is the monomial symmetric function, $h_\mu$ is the
complete symmetric function, and $\delta$ is the usual Kronecker symbol.
Moreover, if $\omega $ is the homomorphism  defined by
$\omega e_i=h_i$ and $\omega h_i=e_i$, where $e_i$ is the elementary
symmetric function, then for any symmetric functions $f$ and $g$ we
have
\begin{equation}
\label{e-omega} \langle f,g\rangle =\langle \omega f, \omega g
\rangle.
\end{equation}

Associate the function
\begin{equation}\label{e-SJ}
S_J=\sum_{\DES w =J} x_{w_1}x_{w_2}\cdots x_{w_n}
\end{equation}
with each subset $J\subseteq [n-1]$, where the sum ranges over all
words on positive integers with descent set $J$. We claim that $S_J$
is a symmetric function whose shape is a border strip (see
\cite[p. 345]{EC2}). In particular, $S_{[n-1]}$ is equal to $e_n$,
the elementary symmetric function of order $n$.
%
On the other hand,
every partition $\lambda$ of $n$ has an associate symmetric function
$L_\lambda$ related to a Lie representation. The definition of
$L_\lambda$ is omitted here (see \cite{gessel-reutenauer});
just remember that the symmetric function
corresponding to derangements of order $n$ is given by
\begin{align}
\mathcal{D}_n&=\sum_{\lambda} L_\lambda =\sum_{j=0}^n (-1)^j e_j
h_1^{n-j}, \label{e-D-0}
\end{align}
where the sum ranges over all partitions $\lambda$
having no part equal to 1
\cite[Theorem 8.1]{gessel-reutenauer}.
We need the following result
from~\cite{gessel-reutenauer} for our enumeration.

\begin{prop}[Gessel-Reutenauer]
\label{t-gessel-reutenauer} The number of permutations having
descent set~$J$ and cycle structure $\lambda$ is equal to the scalar
product of the symmetric functions $S_J$ and~$L_\lambda$.
\end{prop}
\begin{proof}[Proof of Theorem \ref{t-max-fix}]

For each fixed integer sequence $(a_1,a_2,\dots,a_k)$ let
$s_i=a_1+a_2+\cdots+a_i$ for $i=1,\dots,k$ and $\ell=s_k$.
Then
$M(a_1,a_2,\dots,a_k)$ is the number of {\it derangements}
$\pi\in\Sym_\ell$
such that $s_i$ with $i=1,2,\dots,k-1$ may or
may not be a descent of $\pi$, and such that all the other numbers in
$[\ell-1]$ are descents of $\pi$.
There is then a set $T$ of $2^{k-1}$ descent sets $J$ to consider,
depending on whether each $s_i$ is a
descent or not (for $i=1,\dots,k-1$).
By Proposition
\ref{t-gessel-reutenauer} and linearity we have
\begin{equation}
\label{e-D-scalar} M(a_1,a_2,\dots,a_k)= \langle \sum_{J\in T} S_J,
\mathcal{D}_{\ell} \rangle.
\end{equation}
From \eqref{e-SJ} it follows that
$$
\sum_{J\in T} S_J
=\sum_{\DES w\in T }x_{w_1}x_{w_2}\cdots x_{w_n}
=\sum_{[\ell-1]\setminus\{s_1, s_2, \ldots, s_{k-1}\}
\subseteq \DES w }x_{w_1}x_{w_2}\cdots x_{w_n}.
$$
Each word $w$ occurring in the latter sum is the juxtaposition
product $w=u^{(1)} u^{(2)} \cdots u^{(k)}$, where each $u^{(i)}$ is
a decreasing word of length $a_i$ ($i=1,2,\ldots, k$). Hence $
\sum_{J\in T} S_J =e_{a_1}e_{a_2}\cdots e_{a_k}$.
In \eqref{e-D-scalar} replace $\sum_{J\in T} S_J$ by
$e_{a_1}e_{a_2}\cdots e_{a_k}$ and $\mathcal{D}_{\ell}$ by the second
expression in \eqref{e-D-0}. We obtain
$$
M(a_1,a_2,\dots,a_k)= \langle e_{a_1}e_{a_2}\cdots e_{a_k},
\sum_{j=0}^n (-1)^j e_j
h_1^{n-j} \rangle.
$$
The image under  $\omega$ yields
\begin{align*}
M(a_1,a_2,\dots,a_k) &=\langle
\omega e_{a_1}\cdots e_{a_k}, \omega \sum_{j=0}^{\ell} (-1)^j e_j
h_1^{\ell-j} \rangle \\
&= \langle h_{a_1}\cdots h_{a_k},
\sum_{j=0}^{\ell} (-1)^j h_j e_1^{\ell-j}
\rangle.
\end{align*}
Notice that
$\sum_{j=0}^{\ell} (-1)^j h_j e_1^{\ell-j} $
is the coefficient of $u^\ell$ in
$$
\big(\sum_{j}h_j(-u)^j \big)\big(\sum_i e_1^{i}u^i \big) =
\frac{1}{(1+x_1u)(1+x_2u)\cdots (1+x_ku)(1-(x_1+x_2+\cdots +x_k)u)}.
$$
It follows from \eqref{e-scalar} that $M(a_1,\dots,a_k)$ is  the coefficient of
$x_1^{a_1}\cdots x_k^{a_k}u^\ell$  in the
expansion of the
above fraction.
\end{proof}


\section{Permutations with 0 or 1 fixed points}\label{sec-th2} 

Our objective in this section is to prove Theorem \ref{exc-gen}. We
will establish a chain of equivalent or stronger statements, leading
to the final easy one.
Further notations are needed.
Let
$w=w_1w_2\cdots w_n$ be a word on the letters $1,2,\dots,m$, each letter
appearing at least once.
The set-statistic $\IDES w$ is defined to be the set of
all $i$ such that the rightmost $i$ appears to the right of the
rightmost $i+1$ in $w$. Note that if $\pi$ is a permutation on
$[n]$, then $\IDES \pi = \DES \pi^{-1}$. For every proper subset $J$
of $[n-1]$ let $\Sym_n^J$ be the set of permutations $\sigma \in
\Sym_n$ with $\IDES\sigma=J$. Note the difference with the notation
of $D^J_k(n)$ for $k=0,1$. We will see that it is easier to deal
with $\IDES$ than with $\DES$ directly.

A word $w=w_1w_2\cdots w_n$ is said to be a {\it desarrangement} if
$w_1>w_2>\cdots >w_{2k}$ and $w_{2k}\leq w_{2k+1}$ for some~$k\ge
1$. By convention, $w_{n+1}=\infty$. We may also say that the {\it
leftmost trough} of~$w$ occurs at an {\it even} position
\cite{FH07}. This notion was introduced, for permutations, by
D\'esarm\'enien \cite{De84} and elegantly used in a subsequent paper
\cite{DW88}. A further refinement is due to Gessel \cite{Ge91}. A
desarrangement $w=w_1w_2\cdots w_n$ is called a {\it hook}, if $n\ge
2$ and $w_{1}>w_{2}\leq w_{3}\leq \cdots \leq w_{n}$. Every nonempty
word~$w$ on the letters $1,2,3,\dots$ can be written uniquely as a
product $uh_1h_2\cdots h_k$, where $u$ is a  {\it weakly increasing}
word (possibly empty) and each~$h_i$ is a hook. This factorization
is called the {\it hook-factorization} of~$w$ \cite{FH07}.
For permutations it was already introduced  by Gessel
\cite{Ge91}.
For instance, the hook-factorization of the following word is
indicated by vertical bars:
$$
w=\mid 1\,2\,4\,5\mid6\,4\,5\,6\mid 4\,1\,3\mid 6\,5\mid 5\,4\mid
6\,1\,1\,4\mid5\,1\,1\mid.
$$

Let $uh_1h_2\cdots
h_k$ be the hook factorization of the word $w$. The statistic
$\pix w$ is defined to be the length of $u$, and the statistic
$\lec w$ is defined, in terms of inversion statistics ``$\inv$", by
the sum 
\cite{FH07}
$$
\lec w:=\sum_{i=1}^k \inv(h_i).
$$
In the previous example, $\pix w=|1245|=4$ and $\lec
w=\inv(6456)+\inv(413)+\inv(65)+\inv(54)+\inv(6114)+\inv(511)=2+2+1+1+3+2=11$.
\medskip

For each permutation $\sigma$ let $\iexc\sigma=\exc \sigma^{-1}$.
The next proposition was proved in Foata and Han \cite{FH07}.
\begin{prop}\label{FHequi}
The two three-variable statistics $(\iexc, \fix, \IDES)$ and $(\lec,
\pix, \IDES)$ are equi-distributed on the symmetric group $\Sym_n$.
\end{prop}
Let $\Desar_0^J(n)$ denote the set of all desarrangements in
$\Sym_n^J$, and $\Desar_1^J(n)$ the set of all permutations in
$\Sym_n^J$ with exactly one pixed point. Since the map $\sigma\to
\sigma^{-1}$ preserves the number of fixed points, Theorem
\ref{exc-gen} is equivalent to asserting that
$$\sum_{\genfrac{}{}{0pt}{3}{\sigma\in D_0(n)}{\IDES(\sigma)=J}}-
\sum_{\genfrac{}{}{0pt}{3}{\sigma\in D_1(n)}{\IDES(\sigma)=J}}=(s-1)Q_n^J(s). $$
Then by Proposition \ref{FHequi} this is equivalent to the following
Theorem \ref{exc-gen'}.

\addtocounter{thm}{-1}
\renewcommand\thethm{\ref{exc-gen}$^a$}
\begin{thm}\label{exc-gen'}
We have
$$
\sum_{\sigma\in\Desar_0^J(n)} s^{\lec\sigma}
-\sum_{\sigma\in\Desar_1^J(n)} s^{\lec\sigma} =(s-1) Q_n^J(s),
$$
where $Q_n^J(s)$ is a polynomial with positive integral coefficients.
\end{thm}
\renewcommand\thethm{\arabic{thm}}
The following lemma enables us to prove a stronger result.
\begin{lem}\label{l-cases}
Let $w=w_1w_2\cdots w_n$ be a desarrangement such that $\IDES
w\not=\{1,2,\ldots, n-1\}$ and let $w'=w_nw_1w_2\cdots w_{n-1}$. Then,
either $\lec w'=\lec w$, or $\lec w'=\lec w -1$.
\end{lem}
\begin{proof}
Several cases are to be considered.
Say that $w$ belongs to type $A$ if
$\lec(w')=\lec(w)$, and say that $w$ belongs to type $B$ if
$\lec(w')=\lec(w)-1$.

Since $w$ is a desarrangement, we may assume $w_1>w_2>\cdots
>w_{2k}\le w_{2k+1}$ for some $k$. It follows that $w'$ has one pixed point.
Let $h_1\cdots h_k$ be the hook-factorization of $w$. Then the
hook-factorization of $w'$ must have the form $w_n| h_1'\cdots
h_\ell'$. Thus, when computing $\lec(w')$, we can simply omit $w_n$.
This fact will be used when checking the various cases. The reader
is invited to look at Figures 1--3, where the letters $b,c,x,y,z$
play a critical role.

\begin{enumerate}
\item If the rightmost hook $h_k$ has at least three elements, as shown in Figure
\ref{f-case1}, then $b\le c$ belongs to type $A$ and $b> c$ belongs
to type $B$. This is because the only possible change for ``lec"
must come from an inversion containing $c$. Furthermore, $(b,c)$ forms an
inversion for type $B$ and does not form an inversion for type $A$.

{
\def\dy{0.7}
\def\dysize{2.8mm} 
\vskip -\dysize
\def\dot{\circle*{.2}}
\def\putdot(#1,#2){\put(#1,#2){\dot}}
\def\R(#1,#2){\put(#1,#2){\line(1,1){1}}\putdot(#1,#2)}
\def\D(#1,#2){\put(#1,#2){\line(1,-1){1}}\putdot(#1,#2)}
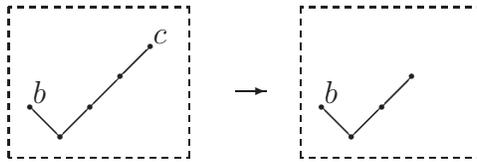
\begin{figure}[hbt]
\begin{center}
\setlength{\unitlength}{4mm}
\begin{picture}(6,5)
\D(0,1) \R(1,0) \R(2,1) \R(3,2)
\putdot(4,3)
\put(0.1,1.1){$b$}
\put(4.1,3.1){$c$}
\put(-0.7,-\dy){\dashbox{.2}(6,5){}}
\end{picture}
\begin{picture}(3,5)
\put(0.5,1.5){\vector(1,0){1}}
\end{picture}
\begin{picture}(6,5)
\D(0,1)\R(1,0)\R(2,1)
\putdot(3,2)
\put(0.1,1.1){$b$}
\put(-0.7,-\dy){\dashbox{.2}(6,5){}}
\end{picture}
\end{center}
\vskip \dysize
\caption{\label{f-case1}Transformation for case $1$.}
\end{figure}
}

\item Suppose the rightmost hook $h_k$ has two elements $b>c$.
 \begin{enumerate}
   \item If there is a hook $xy$ followed by several decreasing hooks of length $2$
with $y\le z$, as shown in Figure \ref{f-case2a}, then $x\le z$
belongs to type $B$ and $x>z$ belongs to type $A$.

{
\def\dy{0.0}
\def\dysize{0mm} 
\vskip -\dysize
\def\dot{\circle*{.2}}
\def\putdot(#1,#2){\put(#1,#2){\dot}}
\def\R(#1,#2){\put(#1,#2){\line(1,1){1}}\putdot(#1,#2)}
\def\D(#1,#2){\put(#1,#2){\line(1,-1){1}}\putdot(#1,#2)}
\begin{figure}[hbt]
\begin{center}
\setlength{\unitlength}{4mm}
\begin{picture}(8,7)
\D(0,4)\putdot(1,3) \put(-0.5,3.3){$x$}\put(0.9,2.3){$y$}
\D(2,6)\putdot(3,5) \put(2.2,6.1){$z$}
\D(4,4)\putdot(5,3)
\D(6,2)\putdot(7,1) \put(6.1,2.1){$b$}\put(7.1,1.1){$c$}
\put(-0.8,-\dy){\dashbox{.2}(9,7){}}
\end{picture}
\begin{picture}(4,7)
\put(1.5,3.0){\vector(1,0){1}}
\end{picture}
\begin{picture}(8,7)
\D(0,4)\putdot(1,3) \put(-0.5,3.3){$x$}\put(0.9,2.3){$y$}
\put(1,3){\line(1,3){1}}
\putdot(2,6)
\D(3,5)\putdot(4,4)
\D(5,3)\putdot(6,2)
\put(2.2,6.1){$z$}
\put(6.1,2.1){$b$}
\put(-0.8,-\dy){\dashbox{.2}(9,7){}}
\end{picture}
\end{center}
\vskip \dysize
\caption{\label{f-case2a}Transformation for case $2a$.}
\end{figure}
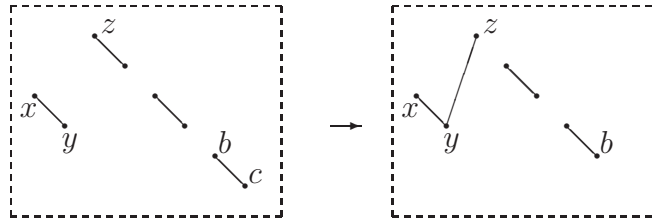
}

\item If
there is a hook of length at least $3$, followed by several
decreasing hooks of length $2$, then (see Figure
\ref{f-case2b})
  \begin{enumerate}
   \item  $x>y$ belongs to type $B$ and $ x\le y$
   belongs to type $A$ in case $y>z$;

   \item $x\le z$ belongs to type $B$ and $x>z$
   belongs to type $A$ in case $y\le z$.
  \end{enumerate}

{
\def\dot{\circle*{.2}}
\def\putdot(#1,#2){\put(#1,#2){\dot}}
\def\R(#1,#2){\put(#1,#2){\line(1,1){1}}\putdot(#1,#2)}
\def\D(#1,#2){\put(#1,#2){\line(1,-1){1}}\putdot(#1,#2)}
\def\boxhsize{7cm}
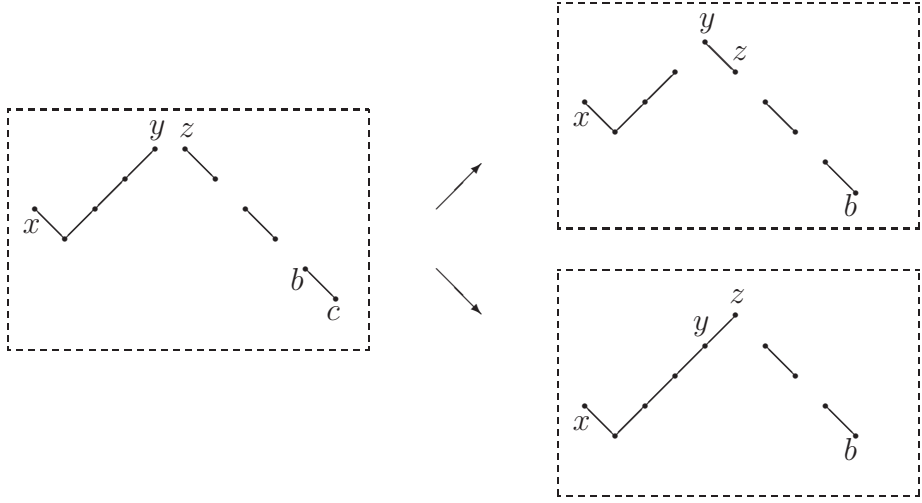
\begin{figure}[hbt]
\begin{center}
\setlength{\unitlength}{4mm}
\hbox{
\hskip 5mm
$\vcenter{
\hsize=\boxhsize
\begin{picture}(11,7)
\D(0,4)\R(1,3)\R(2,4)\R(3,5)\putdot(4,6)
\D(5,6)\putdot(6,5)
\D(7,4)\putdot(8,3)
\D(9,2)\putdot(10,1)
\put(-.4,3.2){$x$}
\put(3.8,6.5){$y$}
\put(4.8,6.4){$z$}
\put(8.5,1.3){$b$}
\put(9.7,0.3){$c$}
\put(-0.9,-0.7){\dashbox{.2}(12,8){}}
\end{picture}
}
$
\hskip -13mm
$\vcenter{
\hsize=2cm
\begin{picture}(4,7)
\put(1.5,4.0){\vector(1,1){1.5}}
\put(1.5,2.0){\vector(1,-1){1.5}}
\end{picture}
}$
\hskip 5mm
$\vcenter{
\hsize=\boxhsize
\hbox{
\begin{picture}(11,10)
\D(0,4)\R(1,3)\R(2,4)\putdot(3,5)
\D(4,6)\putdot(5,5)
\D(6,4)\putdot(7,3)
\D(8,2)\putdot(9,1)
\put(-.4,3.2){$x$}
\put(3.8,6.5){$y$}
\put(4.9,5.4){$z$}
\put(8.6,0.2){$b$}
\put(-0.9,-0.2){\dashbox{.2}(12,7.5){}}
\end{picture}
}
%
%
\hbox{
\begin{picture}(11,10)
\D(0,4)\R(1,3)\R(2,4)\R(3,5)\R(4,6)\putdot(5,7)
\D(6,6)\putdot(7,5)
\D(8,4)\putdot(9,3)
\put(-.4,3.2){$x$}
\put(3.6,6.5){$y$}
\put(4.8,7.4){$z$}
\put(8.6,2.2){$b$}
\put(-0.9,1){\dashbox{.2}(12,7.5){}}
\end{picture}
}
}$
}
\end{center}
\caption{\label{f-case2b}Transformations for case $2b$.}
\end{figure}
}

 \end{enumerate}
\end{enumerate}
This achieves the proof of the lemma.
\end{proof}

With
the notations of Lemma \ref{l-cases} we say that a
desarrangement $w$ is in class $A_0$ if $\lec w'=\lec w$ and in
class $B_0$ if $\lec w'=\lec w -1$. A word $w=w_1w_2w_3\cdots w_n$
is said to be in class $A_1$ (resp. in class $B_1$) if the word
$w_2w_3\cdots w_n w_1$ is in class $A_0$ (resp. in class $B_0$).
Notice that a word in class $A_1$ or $B_1$ has exactly one pixed
point. Then, Theorem \ref{exc-gen'} is a consequence of the
following theorem.

\goodbreak
\addtocounter{thm}{-1}
\renewcommand\thethm{\ref{exc-gen}$^b$}
\begin{thm}\label{exc-gen-AB}
We have
$$
\sum_{\sigma\in\Sym_n^J\cap A_0} s^{\lec\sigma}
 =\sum_{\sigma\in\Sym_n^J\cap A_1} s^{\lec\sigma}
\hbox{\quad and\quad}
\sum_{\sigma\in\Sym_n^J\cap B_0} s^{\lec\sigma}
 =s\sum_{\sigma\in\Sym_n^J\cap B_1} s^{\lec\sigma}.
$$
\end{thm}
\renewcommand\thethm{\arabic{thm}}
\goodbreak

Let $\Sym_n^{\subseteq J}$ be the set of all permutations $\sigma$
of order $n$ such that $\IDES\sigma\subseteq J$. By the
inclusion-exclusion principle,
Theorem \ref{exc-gen-AB} is
equivalent to the following theorem.
\addtocounter{thm}{-1}
\renewcommand\thethm{\ref{exc-gen}$^c$}
\begin{thm}\label{exc-gen-AB'}
We have
$$
\sum_{\sigma\in\Sym_n^{\subseteq J}\cap A_0} s^{\lec\sigma}
=\sum_{\sigma\in\Sym_n^{\subseteq J}\cap A_1} s^{\lec\sigma}
\hbox{\quad and\quad}
\sum_{\sigma\in\Sym_n^{\subseteq J}\cap B_0} s^{\lec\sigma}
=s\sum_{\sigma\in\Sym_n^{\subseteq J}\cap B_1} s^{\lec\sigma}.
$$
\end{thm}
\renewcommand\thethm{\arabic{thm}}

If $J=\{j_1, j_2, \ldots, j_{r-1}\}\subseteq [n-1]$, define a
composition ${\bf m}=(m_1, m_2, \ldots, m_r)$ by $m_1=j_1,
m_2=j_2-j_1, \ldots, m_{r-1}= j_{r-1}-j_{r-2}, m_r=n-j_{r-1}$. Let
$R({\bf m})$ be the set of all rearrangements of $1^{m_1} 2^{m_2}
\cdots r^{m_r}$. We construct a bijection $\phi$ from $R({\bf m})$
to $\Sym_n^{\subseteq J}$ by means of the classical {\it
standardization} of words. Let $w\in R({\bf m})$ be a word. From
left to right label the letters $1$ in $w$ by $1, 2, \dots,m_1$,
then label the letters $2$ in $w$ by $m_1+1,m_1+2,\dots,m_1+m_2$, and so
on. Then the standardization of $w$, denoted by $\sigma=\phi(w)$, is
the permutation obtained by reading those labels from left to right.
It is easy to see that $\phi$ is reversible and
$\IDES\sigma\subseteq J$ if and only if $w\in R(\mathbf{m})$ (see
\cite{DW93, FH07}).
Moreover, the permutation $\sigma$ and the word $w$ have the same
hook-factorization {\it type}. This means that if $ah_1h_2\ldots
h_s$ (resp. $bp_1p_2\ldots p_k$) is the hook-factorization of
$\sigma$ (resp. hook-factorization of $w$), then $k=s$ and
$|a|=|b|$. For each $1\leq i \leq k$ we have $|h_i|=|p_i|$ and
$\inv(h_i)=\inv(p_i)$. Hence $ \lec w=\lec\sigma$ and $ \pix
w=\pix\sigma$. Furthermore, $\sigma$ is in class $A_0$, $A_1$, $B_0$
or $B_1$ if and only if $w$ is in the same class. Theorem
\ref{exc-gen-AB'} is equivalent to the next theorem,
whose proof follows from the definition of the classes $A_0$,
$A_1$, $B_0$, $B_1$ and Lemma \ref{l-cases}.

\addtocounter{thm}{-1}
\renewcommand\thethm{\ref{exc-gen}$^d$}
\begin{thm}\label{exc-gen-AB''}
We have
$$
\sum_{\sigma\in R({\bf m})\cap A_0} s^{\lec\sigma}
 =\sum_{\sigma\in  R({\bf m})\cap A_1} s^{\lec\sigma}
\hbox{\quad and\quad}
\sum_{\sigma\in  R({\bf m})\cap B_0} s^{\lec\sigma}
 =s\sum_{\sigma\in  R({\bf m})\cap B_1} s^{\lec\sigma}.
$$
\end{thm}
\renewcommand\thethm{\arabic{thm}}

The following variation of Theorem \ref{exc-gen} follows from
Theorem \ref{exc-gen-AB}, but cannot be derived from Theorem
\ref{exc-gen} directly.
\begin{thm}\label{exc-gen-var}
We have
$$
\sum_{\sigma\in\Der_0^J(n)} s^{\iexc\sigma}
-\sum_{\sigma\in\Der_1^J(n)} s^{\iexc\sigma} =(s-1) Q_n^J(s)
$$
for some polynomial $Q_n^J(s)$ with positive integral coefficients.
\end{thm}

\section{Further remarks} 
A combinatorial proof of Corollary \ref{t-2-gen} can be made by using the
methods developed in the preceding section.
However this
proof does not need
the concept of ``hook" and the statistic ``$\lec$". We only list the
equivalent statements, leaving the details to the reader.

\begin{thm}\label{t-S2}
Let $J$ be a proper subset of $[n-1]$. The following statements are
equivalent to Corollary \ref{t-2-gen}:
\begin{enumerate}
\item The number of derangements
in $\Sym_n^J$ is equal to the number of permutations in $\Sym_n^J$
with exactly one fixed point.

\item The number of desarrangements in $\Sym_n^J$ is equal to the number
of permutations in $\Sym_n^J$ with exactly one pixed point.

\item The number of desarrangements in $\Sym_n^{\subseteq J}$ is equal to
the number of permutations in $\Sym_n^{\subseteq J}$ with exactly one
pixed point.

\item The number of desarrangements in $R{(\bf m)}$ is equal to the number
of words in  $R{(\bf m)}$ with exactly one pixed point.
\end{enumerate}
 \end{thm}
We remark that the equivalence of $(1)$ and $(2)$ also follows from
a result of D\'esarm\'enien and Wachs \cite{DW88, DW93}: the two
bi-variable statistics (fix, IDES) and (pix, IDES) are
equi-distributed on the symmetric group $\Sym_n$.

\medskip
The statistics ``des" and ``maj" are determined by ``DES": $\des
\pi =\# \DES \pi$ and $ \maj \pi = \sum_{i\in \DES \pi} i$ for $\pi
\in \Sym_n$.
 By using
Theorem \ref{exc-gen}  for each proper subset $J$ of $[n-1]$ and by
checking the case $J=[n-1]$ directly, we have
the following result.
\begin{thm}\label{exc-sum}
There is a polynomial $Q_n(s,t,q)$ with positive integral
coefficients such that
$$
\sum_{\sigma\in\Der_0(n)} s^{\exc\sigma} t^{\des\sigma}
q^{\maj\sigma} -\sum_{\sigma\in\Der_1(n)} s^{\exc\sigma}
t^{\des\sigma} q^{\maj\sigma} =(s-1) Q_n(s,t,q) + r_n(s,t,q)
$$
where $r_{2k}(s,t,q)=s^k t^{2k-1} q^{k(2k-1)}$ for $k\geq 1$ and
$r_{2k+1}(s,t,q)=-s^k t^{2k} q^{k(2k+1)}$ for $k\geq 0$.
\end{thm}

A related result is the following, where we use the standard
notation for $q$-series:
$$(z;q)_m= (1-z)(1-zq)\cdots
(1-zq^{m-1}).$$
\begin{prop}[\cite{FH07}, Theorem 1.1] \label{FH-gf}
Let $(A_{n}(s,t,q,Y))_{n\ge 0}$
 be the sequence of polynomials in
four variables, whose factorial generating function is given by
$$
\sum_{r\ge 0}t^r\frac{(1-sq)\,(u;q)_r\,(usq;q)_r}
{((u;q)_r-sq(usq;q)_r)(uY;q)_{r+1}}\!=\!\sum_{n\ge 0} A_n(s,t,q,Y)
\frac{u^n}{ (t;q)_{n+1}}.
$$
Then $A_{n}(s,t,q,Y)$ is the generating polynomial for~$\Sym_n$
according to the four-variable statistic $(\exc,\des,\maj,\fix)$. In
other words,
$$A_n(s,t,q,Y)=\sum_{\sigma\in\Sym_n}
s^{\exc\sigma}t^{\des\sigma}q^{\maj\sigma} Y^{\fix\sigma}.
$$
\end{prop}

Since $\sum_{\sigma\in\Der_0(n)} s^{\exc\sigma} t^{\des\sigma}
q^{\maj\sigma}$ is simply $A_n(s,t,q,0)$ and
$\sum_{\sigma\in\Der_1(n)} s^{\exc\sigma} t^{\des\sigma}
q^{\maj\sigma}$ is equal to the coefficient of $Y$ in
$A_n(s,t,q,Y)$, Theorem \ref{exc-sum} and Proposition \ref{FH-gf} imply
the following theorem.

\begin{thm}\label{exc-sum-n}
There is a sequence of polynomials $(Q_n(s,t,q))_{n\ge 0}$ with
positive integral coefficients such that
\begin{multline*}
\sum_{r\ge
0}t^r\left(1-u\frac{1-q^{r+1}}{1-q}\right)\frac{(1-sq)\,(u;q)_r\,(usq;q)_r} {((u;q)_r-sq(usq;q)_r)  } -\frac{1}{1-t}\\
= (s-1) \sum_{n\geq 1}Q_n(s,t,q) \frac{u^n}{(t;q)_{n+1}}  +
r(s,t,q),
\end{multline*}
where
$$
r(s,t,q)=\sum_{k\geq 1} s^k t^{2k-1} q^{k(2k-1)}
\frac{u^{2k}}{(t;q)_{2k+1}} -\sum_{k\geq 0} s^k t^{2k} q^{k(2k+1)}
\frac{u^{2k+1}}{(t;q)_{2k+2}}.
$$
\end{thm}

In the case of $t=1$ and $q=1$ the above theorem yields the following corollary.
\begin{cor}\label{exc-sum-n-q1}
For each $n\geq 0$
let $Q_n(s)$ be the coefficient of $u^n/n!$ in the
Taylor expansion of
$$
H(s)=\frac{u-1}{se^{us}-s^2e^u} - \frac{1}{2s\sqrt{s}}\Bigl(
\frac{e^{u\sqrt{s}}}{\sqrt{s}+1}
+\frac{e^{-u\sqrt{s}}}{\sqrt{s}-1}
\Bigr),
$$
that is
$$
H(s)
=\frac{u^3}{3!}
+(s+3)\frac{u^4}{4!}
+(s^2+17s+4)\frac{u^5}{5!}
+(s^3+46s^2+80s+5)\frac{u^6}{6!}+\cdots+Q_n(s)\frac{u^n}{n!}+\cdots
$$
Then, the coefficients $Q_n(s)$ are polynomials in $s$  with positive
integral coefficients.
\end{cor}

It is easy to show that $Q_{2n-1}(1)=D_{2n-1}/2$ and
$Q_{2n}(1)=(D_{2n}-1)/2$ for $n\ge 2$.
By Formula (6.19) in \cite{FH06IV} we have
$$
Q_n(1)=\sum_{2\leq 2k\leq n-1} k \times n(n-1)(n-2)\cdots (2k+2).
$$
Since $Q_{n}(1)$ counts the number of desarrangements of type $B$,
Corollary \ref{exc-sum-n-q1} implies that the number of
desarrangements of type $A$ equals the number of desarrangements of
type $B$, when excluding the decreasing desarrangement of even
length. It would be interesting to have a direct (analytic) proof of Corollary
\ref{exc-sum-n-q1} which would not use the combinatorial set-up of this paper.

\goodbreak

\vspace{.2cm} \noindent{\bf Acknowledgments.} The authors would like
to thank Dominique Foata for helpful remarks and suggestions. This
work was supported by the 973 Project, the PCSIRT project of the
Ministry of Education, the Ministry of Science and Technology and
the National Science Foundation of China.


\end{document}